\documentclass[12pt]{article}
\usepackage{amsmath}
\usepackage{amssymb}

\newlength{\defbaselineskip}
\newcommand{\setlinespacing}[1]%
           {\setlength{\baselineskip}{#1 \defbaselineskip}}
\newcommand{\doublespacing}{\setlength{\baselineskip}%
                           {2.0 \defbaselineskip}}

\newcommand{\QED}{\hspace*{\fill}\rule{2.5mm}{2.5mm}}

\newenvironment{proof}{\noindent{\bf Proof.\ }}{\QED\\}

\newcommand{\N}{\mathbb N}

\newcommand{\R}{\mathbb R}

\newcommand{\Z}{\mathbb Z}
\def\N{\mathbb N}

\def \proclaim#1{\smallskip\noindent{\bf\ignorespaces#1\unskip.}%
 \bgroup\it\space\ignorespaces}
\def \endproclaim{\par\egroup\smallskip}

\begin{document}

\title{\bf Symmetric polynomials and $l^p$ inequalities for certain intervals of $p$. }
\date{}
\author{Ivo Kleme\v{s} }
\maketitle
{\centerline
{\it
\noindent Department of Mathematics and Statistics,
\noindent 805 Sherbrooke Street West,} }

{\centerline
{\it
\noindent McGill University,
\noindent Montr\'eal, Qu\'ebec,
\noindent H3A 2K6,
\noindent Canada.
}}
\medskip

{\centerline{  Email:  klemes@math.mcgill.ca }}

\bigskip

\bigskip

\noindent {\it Abstract.} We prove some sufficient conditions
implying $l^p$ inequalities of the form $||x||_p \leq ||y||_p$ for vectors
$ x, y \in [0,\infty)^n$ and for $p$ in certain positive real intervals.
Our sufficient conditions are strictly weaker than the usual majorization
relation.
The conditions are expressed in terms
of certain homogeneous symmetric polynomials in the entries of the vectors. These
polynomials include the elementary symmetric polynomials as a special case.
We also give a characterization of the majorization relation
by means of symmetric polynomials.


\vfill
\noindent {\it A.M.S. Mathematics Subject Classifications:} 47A30 (26B25, 52A40).
\bigskip


\bigskip

\noindent  {\it Date:} February 2010. Revised 18 April 2010.

\newpage
%
%


\doublespacing

 \centerline{\bf \S 1. Introduction}
\medskip


Let $x$ and $y$ be given vectors in $\R^n$ having nonnegative entries.
We will investigate sufficient conditions on $x$ and $y$ for
$l^p$ inequalities of the form
$||x||_p \leq ||y||_p$ simultaneously for all $0 \leq p \leq 1$.
Under the additional assumption that $||x||_1 = ||y||_1$, our conditions
also imply
  $||x||_p \geq ||y||_p$ for $1 \leq p \leq 1+r$, where
  $r \geq 1$ is a freely adjustable
integer parameter appearing in the conditions.
As will be seen in Theorem 1, the conditions are expressed using
a finite number of symmetric polynomials in $x$ or $y$ with positive
coefficients, whose degrees are controlled by $r$ in some way. In particular,
the special case $r=1$ of these conditions involves just the
elementary symmetric polynomials. This case is a kind of ``folk theorem".
It has typically been used in order to obtain $l^p$ estimates for
the eigenvalues of some operator $A$, via the determinant of $(I+tA)$
\cite[Ch. 4, p. 211-212, Lemma 11.1]{GK}, \cite[Theorem 4]{MS},
\cite[Theorem 1.2]{k1}.

Such polynomial conditions may be viewed as expressing
certain averaged properties of the $k$th tensor powers
$x^{\otimes k}$ and $y^{\otimes k}$ for various $k$.
As a complement to Theorem 1, we will present in \S3 an almost trivial
characterization of the usual majorization relation $x \succ y$ from the same
point of view, that is by means of certain symmetric polynomials in $x$ or $y$
(Theorem 2). More precisely, we supply a converse to
a previous result by
Proschan and Sethuraman \cite[Theorem 3.J.2, Example 3.J.2.b]{MO}
regarding a class of Schur-concave symmetric polynomials.

A considerable amount of literature exists concerning the
larger set of simultaneous $l^p$ inequalities
given by $||x||_p \leq ||y||_p$ for $-\infty \leq p \leq 1$
and $||x||_p \geq ||y||_p$ for $1 \leq p \leq \infty$. This relation is implied by,
but strictly weaker than $x \succ y$, and has been called
``power majorization" \cite{B1}. It has been studied in the context of
some concrete numerical sequences \cite{B2}, \cite{G}, and also in
quantum information theory, where certain characterizations
have recently been obtained \cite{Kli}, \cite{T}, \cite{A1}, \cite{A2}.
It is interesting that the latter quantum information literature is concerned with relations
of the form $x^{\otimes k} \succ y^{\otimes k}$, and also
$x \otimes z \succ y \otimes z$ for some $z$ (the ``catalyst").
However, the characterizations themselves are more in the spirit of existence proofs,
rather than explicit conditions that can be checked in concrete situations.

Theorems 1 and 2 and their proofs were originally presented by the author
in the 2002 preliminary report \cite{k2} along with a number of related results.
This and some further results were submitted to a journal in February 2007 in the form of
preprint \cite{k3}. Two years later (January 2009)
the journal reported that it had been unable to recruit any referees.
Also, during the latter waiting process the author decided to post
\cite{k3} on arXiv (June 2008).

\bigskip

\centerline{\bf \S 2. The main result. }
\medskip

Let us fix the following notation for the $l^p$ means of a vector $x \in \R^n$:
$||x||_p := \left(\frac{1}{n}\sum_{i=1}^n |x_i|^p \right)^{1/p}\ , 0 \neq p \in \R $
(with the convention that for a negative $p$ we set $||x||_p = 0$ whenever some
entry $x_i=0$, as would be demanded by continuity in $x$),
$||x||_{-\infty} := \min_i |x_i|\ $,
$||x||_0 := \left(\prod_{i=1}^n |x_i| \right)^{1/n}\ $, $||x||_\infty := \max_i |x_i| $,
as demanded by continuity in $p$.

\proclaim {Definition 1 } Let $r \geq 1$ be an integer.
Let $P_r$ be the $r$th degree Taylor polynomial
of $\exp$, that is $P_r(s) = 1 + \frac{s^1}{1!} + \dots + \frac{s^r}{r!}$.
If $x \in \R^n$ and $t \in \R$ let
\begin{equation}
\label{gen1}
f_r(x,t):=
\prod_{i=1}^n \ P_r(x_it) =
\prod_{i=1}^n \left(1 + x_it + \dots + \frac{(x_it)^r}{r!}
\right).
\end{equation}
For each integer $k \geq 1$ define $F_{k,r}(x)$ to be the coefficient of $t^k$ in $f_r(x,t)$.
\endproclaim
\bigskip

\noindent Note that we have not explicitly indicated $n$ in the notation $f_r$ and
$F_{k,r}$, but this should not cause any confusion.
Clearly, the $F_{k,r}$ can be written out
explicitly as
\begin{equation}
\label{Fkr}
F_{k,r}(x_1,\dots,x_n)\  = \ \ \sum_{\sum k_i = k, \ \max k_i \leq r}
\ \ \prod_{i=1}^n \frac{x_i^{k_i}}{k_i!}\ ,
\end{equation}
where it is understood that the $(k_i)_{i=1}^n$ range over
$n$-tuples of nonnegative integers. Equivalently,
$F_{k,r}(x)$ is the sum of those terms in the
expansion of $\frac{1}{k!}(x_1 + \dots + x_n)^k$ in which each
variable $x_i$ has exponent at most $r$.
Clearly $F_{k,1} =: E_k $ is the elementary symmetric polynomial
of degree $k$ and $F_{k,r} = (E_1)^k/k! \ $
whenever $k \leq r$. Also,
$F_{r+1,r}(x) = (E_1(x)^{r+1} - \sum_i x_i^{r+1})/{(r+1)}!$,
$ F_{nr,r}(x) = (E_n(x))^r/(r!)^n$, and $F_{k,r}(x)=0$  when $k > nr$.
Our main result is the following.
\bigskip

\proclaim { Theorem 1}
Let $x,y \in [0,\infty)^n$ and fix an integer $r \geq 1. $
If
\begin{equation}
\label{F0}
 F_{k,r}(x) \leq  F_{k,r}(y)
 \end{equation}
for all integers  $k$
in the interval $r \leq k \leq nr,$ then
\begin{equation}
\label{p0}
 ||x||_p \leq ||y||_p \ \ \ {\it whenever} \ \ \ 0 \leq p \leq 1.
\end{equation}
If also $ {\displaystyle \sum_{i=1}^n x_i  = \sum_{i=1}^n y_i  }$, then
\begin{equation}
\label{p1}
 ||x||_p \geq ||y||_p \ \ \ {\it whenever} \ \ \ 1 \leq p \leq r+1.
\end{equation}
\endproclaim
%
\bigskip

\begin{proof} Fix the integer $r\geq1.$
Observe that $\log(1 + s + \dots + \frac{s^r}{r!})$ is $O(s)$ when
$s\to 0^{+}$ and $O(\log s)$ when $s \to +\infty.$
Thus, the  integrals (Mellin transforms)
$$I_r(p)  := \int_0^\infty \log(1 + s + \dots + \frac{s^r}{r!})
\ s^{-p} \  \frac{ds}{s}$$
are finite (and positive) for all $p$ in the interval $0<p<1.$ Replacing $s$
by $at$ for any positive $a$ gives the identity
\begin{equation}
\label{id1}
\frac{1}{I_r(p)}\int_0^\infty \log(1 + at + \dots + \frac{(at)^r}{r!})
\ t^{-p} \  \frac{dt}{t} \ = \ a^p \ \ \  \ (a \geq 0,\ 0<p<1).
\end{equation}
 Now let
$x,y \in [0,\infty)^n$ and
 $ F_{k,r}(x) \leq  F_{k,r}(y) $
for all integers  $k$
in the interval $r \leq k \leq nr$. Note that in the case $r=k$ we have
$F_{r,r}(x) = (E_1(x))^r/r! = (n||x||_1)^r/r!$. Hence $||x||_1 \leq ||y||_1$.
Also, $F_{k,r} = (E_1)^k/k! \ $ for $1\leq k \leq r$. Thus
in fact $ F_{k,r}(x) \leq  F_{k,r}(y) $
for all integers  $k$
in the interval $1 \leq k \leq nr$, i.e. for all coefficients of $t^k$ in the
generating functions $f_r(x,t)$ and $ f_r(y,t)$ (see Definition 1).
Hence
$$1\leq f_r(x,t) \leq f_r(y,t) \ , \ \ \ \forall \ \ t \geq 0 .$$
Taking logarithms of the $f_r$ and integrating with respect to $t^{-p} \  \frac{dt}{t}
\frac{1}{I_r(p)}$ gives, by identity (\ref{id1}),
$$ \sum_{i=1}^n x_i^p \leq \sum_{i=1}^n y_i^p  \  \ \ \ \ ( 0<p<1). $$
Normalizing both sides
we obtain the first case of the theorem, since the inequalities
$||x||_p \leq ||y||_p$ extend to the endpoint case $p=0$ by continuity in $p$.
Next, if in addition $\sum_i x_i  = \sum_i y_i \ ,$ then
$\sum_i x_it  = \sum_i y_it \ $ for all $t\geq 0$. Subtracting from this
the inequality $ \log f_r(x,t) \leq \log f_r(y,t)$, one obtains
$$
 \sum_i \left( x_it - \log(1 + x_it + \dots + \frac{(x_it)^r}{r!}) \right) \
$$
\begin{equation} \label{sum}
\geq \  \sum_i \left( y_it - \log(1 + y_it + \dots + \frac{(y_it)^r}{r!}) \right).
\end{equation}
Consider the function $\delta_r(s):= s-\log(1 + s + \dots + \frac{s^r}{r!})$
for $s\geq 0.$
We have $\delta_r(s) \geq s- \log(e^s) =0$ for $s\geq 0.$ When $s\to +\infty,$ we
have $\delta_r(s) = O(s) + O(\log (s^r) ) = O(s).$ When $s \to 0^{+}$ we have
$\delta_r(s) = s - \log\left(e^s-O(s^{r+1})\right) =
s-\log\left(e^s(1-e^{-s}O(s^{r+1})\right) =
s-\log(e^s)-\log\left(1-e^{-s}O(s^{r+1})\right) = O(e^{-s}O(s^{r+1}))
=O(s^{r+1}).$ It follows that the  integrals
$$J_r(p)  := \int_0^\infty
\left(s-\log(1 + s + \dots + \frac{s^r}{r!})\right)
\ s^{-p} \  \frac{ds}{s}$$
are finite (and positive) for all $p$ in the interval $1<p<r+1.$
Replacing $s$
by $at$ gives the  new identity
\begin{equation}
\label{id2}
\frac{1}{J_r(p)}\int_0^\infty
 \left(at-\log(1 + at + \dots + \frac{(at)^r}{r!})\right)
\ t^{-p} \  \frac{dt}{t} \ = \ a^p  \ ,
\end{equation}
for $a\geq0,\ 1<p<r+1$.
Thus, when $1<p<r+1$ we may integrate (\ref{sum}) with respect to
$t^{-p} \  \frac{dt}{t}
\frac{1}{J_r(p)}$ and use (\ref{id2}) to obtain
$$ \sum_i x_i^p \geq \sum_i y_i^p  \ , \ \ \ \ ( 1<p<r+1). $$
By continuity in $p,$ we obtain $||x||_p \geq ||y||_p$ for
$1\leq p \leq r+1.$
\end{proof}

\noindent {\bf Remarks on Theorem 1:}

{\bf (a).} The case $r=1$ of Theorem 1 employs only the elementary symmetric
polynomials $E_k =F_{k,1}$ and is relatively well known, as mentioned
in the introduction. We illustrate the cases $r= 1, 2$
in an example following these remarks.

{\bf (b).} One can ask some natural questions regarding the sharpness of
various aspects of Theorem 1, but we will not go into the details within
the space of the present paper. Let us mention only the following without proof
(some of these remarks are discussed further in \cite{k3}):
(i) In the conclusions (\ref{p0}) and (\ref{p1}), the intervals of $p$
cannot be enlarged at either end, at least when $n \geq 3$. In particular,
one cannot make any general conclusion in the range $p<0$.
(ii) The converse of Theorem 1 does not hold in general, in the sense that
(\ref{p0}) and (\ref{p1}) do not imply the hypotheses (\ref{F0}), when $n \geq 4$.
There is a strong converse when $n=3$ and $\sum x_i = \sum y_i$ : Then the two end point
$l^p$ inequalities $||x||_0 \leq ||y||_0$  and  $||x||_{r+1} \geq ||y||_{r+1}$
imply all of the hypotheses (\ref{F0}) concerning the $F_{k,r}$ for a fixed $r$.
(And hence they also imply all the interior cases of $p$ in (\ref{p0}) and (\ref{p1})).
(iii) For general $n$, although there is no converse,
there may be some redundancy in the hypotheses (\ref{F0}).
That is, perhaps some of the $k$'s can be omitted from the
current list $r \leq k \leq nr$. (iv) When $r$ is increased, do the hypotheses
(\ref{F0}) get stronger ? The conclusions suggest that they do. But
on the other hand, for $r_1 < r_2$ the family of functions $\{F_{k,r_1}\}_{k=1}^\infty$
is not simply a subset of the family $\{F_{k,r_2}\}_{k=1}^\infty$ ; one may need
to examine the convex cones spanned by their gradients to answer the question.

{\bf (c.1).} In Theorem 1 the $F_{k,r}$ can be replaced by different
choices of special polynomials as follows. Fix the index $r \geq 1$.
In the proof, only some key properties of the Taylor polynomial
$P_r(s)=1 + \frac{s^1}{1!} + \dots + \frac{s^r}{r!}$ were needed:
We could have replaced $P_r(s)$ by any expression of the form
$$Q_r(s) :=(1 + \frac{s^1}{1!} + \dots + \frac{s^r}{r!}) +
\sum_{j=r+1}^\infty a_{r,j}\frac{s^j}{j!}$$
for any fixed set of constants $0\leq a_{r,j} <1$ having the property
that $\log Q_r(s) \leq K_\epsilon s^\epsilon$ as $s \to \infty$
for any $\epsilon > 0$, i.e. $\log Q_r(s) =O(s^\epsilon)$
for any $\epsilon > 0$. Thus, $Q_r(s)$ should have ``order zero" in the
sense of entire functions; see for example \cite[Ch. 1]{L}.
[Moreover, even with the weaker property that as $s \to \infty$,
$\log Q_r(s) =O(s^\epsilon)$ for a fixed $1 > \epsilon > 0$,
the proof of Theorem 1 still succeeds for the $l^p$ inequalities
in the range $\epsilon \leq p \leq 1$ for (\ref{p0}), and the full range
$1 \leq p \leq r+1$ for (\ref{p1}).]
 We can then use $Q_r(s)$ to define a new generating
function $f_r(x,t) = \prod_{i=1}^n \ Q_r(x_it)$ and
re-define $F_{k,r}(x)$ to be the coefficient of $t^k$ in $f_r(x,t)$.
Theorem 1 then holds as before (of course, the hypothesis
$ r \leq k \leq nr$ should be loosened to include all $ r \leq k  < \infty$).

{\bf (c.2).} The following are some natural examples of Remark (c.1).
 For simplicity we first consider the case $r=1$.
(i) Notice that the inequality
{$\prod_{i=1}^n (1+x_it) \leq \prod_{i=1}^n (1+y_it)$}
would hold if it was known that
$$\big(\prod_{i=1}^n (1+x_it)\big)^M \leq \big(\prod_{i=1}^n (1+y_it)\big)^M$$
for some fixed integer $M \geq 2$. So, we could consider the coefficients
$\widetilde{E}_{k}(x)$ of $t^k$ in the expansion
of $\big(\prod_{i=1}^n (1+x_it)\big)^M$, instead of the usual elementary symmetric
polynomials $E_k(x)$. The weaker
hypothesis $\widetilde{E}_{k}(x) \leq \widetilde{E}_{k}(y) \ \forall k$
would clearly suffice in the $r=1$ case of Theorem 1.
(ii) More generally consider any finite or infinite product
$Q(s):= \prod_{j=0}^\infty (1+c_js)$ with $c_0=1, c_j \geq 0$ and $c_j \to 0$
sufficiently fast to guarantee that $Q(s)$ converges
and $\log Q(s) =O(s^\epsilon)$ for any $\epsilon > 0$
as $s \to \infty$. This is can be seen to be equivalent to the simple
requirement that the sequence $c=\{c_j\}$ belong to $l^\epsilon$ for
every $\epsilon > 0$ \cite[Ch. 1, \S 5]{L}. (For example, $c_j = q^j$
with $0<q<1$.)
Then let $E_{k,c}(x)$ be the coefficient of $t^k$ in $\prod_{i=1}^n Q(x_it)$.
The hypotheses $E_{k,c}(x) \leq E_{k,c}(y)$ for all $k$
would again suffice in the $r=1$ case of Theorem 1.
For the general $r\geq 1$ in Theorem 1, similar modifications of the $F_{k,r}(x)$ can
be constructed by considering products of the form $Q_r(s):= \prod_{j=0}^\infty P_r(c_js)$
in place of the latter $Q(s)$.

{\bf (c.3).}
We note that the discussion in remark (c.2) is equivalent to considering
the finite or infinite ``catalyst" $c= \{c_j\}$ and
comparing various properties of the
two vectors $x \otimes c = \{x_ic_j\}$ and
$y \otimes c = \{y_ic_j\}$.
(See the quantum
information literature mentioned in the introduction for the background.)
Thus, one sees that
$E_{k,c}(x) = E_{k}(x \otimes c)$ and that
each $E_{k,c}(x)$ is in fact a
certain convolution of the sequence $\{E_m(x)\}_{m=1}^k \ $.

\noindent {\bf Example of Theorem 1:}

Motivated by \cite{k1}, we give an example of both
the applicability and inapplicability of Theorem 1.
Suppose that one is interested in comparing the $l^p$ norms of
the eigenvalues $x:= (x_1,\dots, x_4)$ and $y:= (y_1,\dots, y_4)$
respectively of the $4 \times 4$ matrices $X$ and $Y$ defined by
$X= QQ^T$, $Y=RR^T$, where
$$
Q  =  \left[ \begin{array}{cccc}
            1 & 1 & 1 & 1 \\
            0 & 1 & 1 & 0 \\
            0 & 0 & 1 & 0 \\
            0 & 0 & 1 & 1 \\
          \end{array} \right], \ \
R  =  \left[ \begin{array}{cccc}
            1 & 1 & 1 & 1 \\
            0 & 1 & 1 & 0 \\
            0 & 0 & 1 & 0 \\
            0 & 0 & 1 & -1 \\
          \end{array} \right].
$$
A computer plot of $||x||_p$ and $||y||_p$ versus $p$ seems to indicate
that $||x||_p \leq ||y||_p$ for $0 \leq p \leq 1$, and
that $||x||_p \geq ||y||_p$ for $1 \leq p \leq \infty $. Thus, the natural
question is to ask for an ``enlightening" proof or disproof.
More generally, do these inequalities hold whenever $Q$ is a rectangular
(0,1) ``interval matrix" (the 1's occur in some interval in each row)
and $R$ is obtained by arbitrarily changing
signs in the entries of $Q$ ?
We will see that Theorem 1 can be applied to the above example
in the cases $r=1, 2$, but that it does not apply when $r=3$.
Thus the theorem provides a proof of the conjectured inequalities
in the range $0\leq p \leq 2+1=3$, although they appear to
be true for all higher $p$'s as well.

Considering first the range $0 \leq p \leq 2$,
the case $r=1$ of Theorem 1 provides a reasonably nice proof of the
asserted inequalities: One sees that $\sum x_i = \sum y_i =4+2+1+2=9$
and that $E_k(x) \leq E_k(y), \ k=2,3,4$. The latter can be checked by either
directly calculating all coefficients in the two polynomials
$\sum t^k E_k(x)=\det(I+tX)= 1+9t+16t^2+9t^3+t^4$
and $\sum t^k E_k(y)=\det(I+tY)=1+9t+20t^2+9t^3+t^4$, or more efficiently
(see \cite[Section 2]{k1}),
by noting that $X$ is ``totally unimodular" (since it is an interval matrix)
and that $X \equiv Y$ mod 2.

Next, to apply the case $r=2$ of Theorem 1, we need to check whether
or not $F_{k,2}(x) \leq F_{k,2}(y), \ k= 2,3,4,5,6,7,8$. It is not difficult
to express these $F_{k,2}$ as polynomials in the $E_k$ with rational coefficients,
and thus compute their exact values from the above information
(or one could choose to directly compute them as coefficients in the generating
function $\sum_k F_{k,2}(x)t^k = \det(I + Xt + X^2t^2/2!)$).
The results are that

$\{k!F_{k,2}(x)\}_{k=2}^8 = (81,405,1524,4050,7290,5670,2520)$

$\{k!F_{k,2}(y)\}_{k=2}^8 = (81,513,2388,5130,7290,5670,2520)$.

\noindent Since the required inequalities hold, Theorem 1 applies and thus
the proof of the asserted $l^p$ inequalities has been extended to the range $2\leq p \leq 3$.

Finally, attempting to apply the case $r=3$ of Theorem 1, we run into the problem
that $10!F_{10,3}(x)= 1226400 > 1192800 = 10!F_{10,3}(y)$, so that the
hypotheses of Theorem 1 do not hold. (This incidentally also shows that $x$ does not
majorize $y$ in this example, since the $F_{k,r}$ are Schur concave,
as will be discussed in the next section.)
\bigskip

\centerline{\bf \S 3. Comparisons with the majorization relation. }
\medskip

We may put Theorem 1 into a wider context by observing that each of the
functions $F_{k,r}$ is Schur-concave. We will derive this in Example 1 below, but first
we briefly review the relevant
topics concerning the majorization relation $x\succ y$ (also denoted by
$y \prec x$).
A comprehensive treatment may be found in \cite{MO}.

For $x,y \in [0,\infty)^n  $, we write $x\succ y$
(read $x$ ``majorizes" $y$) if
$\sum_{i=1}^n x_i = \sum_{i=1}^n y_i$
and $\sum_{i=1}^k \ x_i^* \ \geq \ \sum_{i=1}^k \ y_i^* $ for
$k=1, \dots , n-1,$ where $x_1^* \geq \dots \geq x_n^*$ denotes the
decreasing rearrangement of the entries $x_i$ of a vector $x$.
The relation $x\succ y$ is equivalent to
the conditions
$\sum_{i=1}^n \varphi(x_i) \geq \sum_{i=1}^n \varphi(y_i)$ for all
convex $\varphi : [0,\infty)\to \R$ and
$\sum_{i=1}^n x_i = \sum_{i=1}^n y_i$.

A symmetric real-valued function $\Phi$ on $[0,\infty)^n $ is called Schur-concave
if $x \succ y \Rightarrow \Phi(x) \leq \Phi(y) $, and Schur-convex
if $x \succ y \Rightarrow \Phi(x) \geq \Phi(y) $. (Hence, $\Phi$ is Schur-convex
if and only if $-\Phi$ is Schur-concave.)

For smooth $\Phi$,
Schur-concavity is equivalent to the well known
Schur-Ostrowski criterion \cite[Theorems 3.A.7, 3.A.8]{MO}: For every pair $i\neq j$,
\begin{equation}
\label{Sch}
\bigg(\frac{\partial \Phi}{\partial x_i} -
\frac{\partial \Phi}{\partial x_j}\bigg)/(x_j -x_i)
\ \geq  \ 0 \  \ \ \ \forall x \in [0,\infty)^n\ \ {\rm with} \ \ x_i\neq x_j.
\end{equation}
This test works in a particularly satisfying way with certain polynomials $\Phi$ where the
quotient in (\ref{Sch}) simplifies to a new polynomial with positive
coefficients. Examples of such nice polynomials $\Phi$ are the elementary
symmetric polynomials $E_k$, more generally all of the $F_{k,r}$, and even more
generally the polynomials $H_S$ in the following result
of Proschan and Sethuraman.
To state the result, let
$I_k = \{ p=(p_1, \dots, p_n) \in \Z^n \ |\ \sum p_i = k \ {\rm and} \
p_i \geq 0 \ \forall i \}.$
  A subset $S \subset I_k$ is said
 to be a ``Schur-concave index set" if its indicator function
 ${\bf 1}_S$ is Schur-concave on $I_k \ $,
 that is if
 $p \succ q  \Rightarrow  {\bf 1}_S(p) \leq {\bf 1}_S(q)$,
 or equivalently, if
 \begin{equation}
\label{S}
 p \in S, \ q \in I_k,  \ p \succ q
\Rightarrow  q \in S.
\end{equation}
\proclaim{Theorem A [Proschan and Sethuraman]}  \cite[Theorem 3.J.2, Example 3.J.2.b]{MO}.
Let $k,n \geq 1$ and let $S \subset I_k$ be a Schur-concave index set.
Define the polynomial $H_S$ by
\begin{equation}
\label{HS}
H_S(x_1,\dots,x_n)\  = \ \ \sum_{p \in S}
\ \ \prod_{i=1}^n \frac{x_i^{p_i}}{p_i!}\ .
\end{equation}
Then for all $x,y \in [0,\infty)^n  \ , $
$ x \succ y
\Rightarrow  H_S(x) \leq H_S(y), $ that is, $H_S$ is Schur-concave
on $[0,\infty)^n  $.
\endproclaim
(As hinted above, one way to prove this theorem is by directly computing  the quotient in (\ref{Sch})
with $\Phi = H_S$ and seeing that it simplifies to a polynomial with positive coefficients.
The interested reader may either try this as an exercise, or refer to \cite{MO} for a proof.)

{\bf Example 1.} Fix integers $k,r \geq 1$. Let $S= \{ p \in \Z^n \ |\ \sum p_i = k, \
p_i \geq 0 \ \forall i , \ {\rm and} \  \max p_i \leq r \}$. It is easy to see that
$S$ is a Schur-concave index set, hence $H_S$ is Schur-concave by Theorem A.
Clearly, $H_S = F_{k,r}$ by (\ref{Fkr}).

Thus, by the Schur-concavity of the $F_{k,r}$ we see
that the majorization relation $x \succ y $
implies $F_{k,r}(x) \leq F_{k,r}(y) $, i.e. the hypotheses of Theorem 1, for all $r \geq 1$.
The converse is false however (when $n \geq 3$). We omit the details but we can make a
remark analogous to (b)(ii) in \S2 above: When $n=3$ it can be
shown that the three norm conditions $||x||_0 \leq ||y||_0$, $||x||_1 = ||y||_1$,
and $||x||_\infty \geq ||y||_\infty$ suffice to
imply $F_{k,r}(x) \leq F_{k,r}(y) $ for all $r \geq 1$. But these
norm conditions do not imply $x \succ y $ (take $x=(4,2-t, 1+t), \
y=(3,3,1)$ for small $t>0$).

{\bf Example 2.} Fix integers $k,r \geq 1$. Let $ S =
\{ \ p \in \Z^n \ |\sum p_i = k, \
p_i \geq 0 \ \forall i ,  \ {\rm and} \ \ p \
\ {\rm has \ at \ least \ } r
{\rm \ nonzero \ entries} \  \}$.
It is again easy to see that
$S$ is a Schur-concave index set, hence $H_S$ is Schur-concave by Theorem A.
Let us introduce the notation $H_S = G_{k,r}$ for these polynomials.
We may think of $G_{k,r}(x)$ as a sum of certain terms in the multinomial
expansion of $\frac{1}{k!}(x_1 + \dots + x_n)^k$, namely those containing
at least $r$ distinct variables $x_i$ as factors.
(If $k<r$ then $G_{k,r}=0$ by the empty sum convention.)

Unlike the $F_{k,r}$, the polynomials $G_{k,r}$ of Example 2 do characterize
the majorization relation $x \succ y $, as will be seen in Theorem 2 below.
As a bonus we also introduce the following closely related symmetric polynomials:
\begin{equation}
\label{Mkr}
M_{k,r}(x):= \sum_{1\leq i_1 < \dots < i_r \leq n}
\ (x_{i_1} + \dots + x_{i_r})^k  \  \ \ \ \ \ \ \ (k, r \in \N).
\end{equation}
(Thus, if $r>n$ we have $M_{k,r}(x)=0 $ by the empty sum convention).
\medskip

\proclaim{Theorem 2}
Let $x,y \in [0,\infty)^n  $ with $ \sum x_i = \sum y_i $.
Then the following three properties are equivalent:
{\rm (a) }$ x \succ y $,
{\rm (b) }$G_{k,r}(x) \leq G_{k,r}(y) $ for all integers $\ k,r \geq 1$,
{\rm (c) }$M_{k,r}(x) \geq M_{k,r}(y) $ for all integers $\ k,r \geq 1$.
\endproclaim
\begin{proof}
Given $x \in [0,\infty)^n$
and a fixed $1\leq r \leq n,$ we may ``compute" the
function $s_r(x):=\sum_{i=1}^r \ x_i^*$
by first noting that it is the maximum of all possible sums of $r$
entries of $x,$ and then computing this maximum by
using integer $k$-norms as $k \to \infty$ :
$$
s_r(x) = \lim_{k\to \infty} \
\left( \sum_{1\leq i_1 < \dots < i_r \leq n}
\ (x_{i_1} + \dots + x_{i_r})^k \ \right)^{\frac{1}{k}} =
\lim_{k\to \infty} \
\left( M_{k,r}(x) \right)^{\frac{1}{k}}.
$$
Thus clearly (c) implies (a). But (a) implies (b) since the $G_{k,r}$
are Schur-concave. To see that (b) implies (c),
it remains to relate the polynomials $M_{k,r}(x)$ to the
$G_{k,r}(x)$.
Consider the polynomials $\overline{G_{k,r}}$ defined by
$$ \overline{G_{k,r}}(x) \ := \
\frac{1}{k!}(x_1 + \dots + x_n)^k - G_{k,r}(x) \ ,
$$
which may be thought of as the sum of all terms in the
expansion of $\frac{1}{k!}(x_1 + \dots + x_n)^k$ containing
{\it less than} $r$ distinct $x_i$ as factors.
To complete the proof, it suffices to show that when
$1 \leq r \leq n$ each $M_{k,r}$
is a linear combination, with positive coefficients, of some of the
$\overline{G_{k,r}}$.
Let
$$\Delta\overline{G_{k,r}}(x)=\overline{G_{k,r+1}}(x) - \overline{G_{k,r}}(x)$$
i.e. the sum of all terms containing
{\it exactly} $r$ distinct $x_i$ as factors.
An expansion of each power in $M_{k,r}$ by the multinomial
theorem gives
$$\frac{1}{k!}M_{k,r} =
{\small \bigg(\begin{array}{c}
n-r  \\
0
 \end{array}\bigg) }
\Delta\overline{G_{k,r}} \ + \
{\small \bigg(\begin{array}{c}
n-r+1  \\
1
 \end{array}\bigg) }
\Delta\overline{G_{k,r-1}}
\ + \dots + \
{\small \bigg(\begin{array}{c}
n-1  \\
r-1
 \end{array}\bigg) }
\Delta\overline{G_{k,1}} \ \ .
$$
Since these binomial coefficients are increasing from left to right,
the result follows after a summation by parts. In fact,
by Pascal's identity we obtain the explicit formula
$$\frac{1}{k!}M_{k,r} =
{\small \bigg(\begin{array}{c}
n-r-1  \\
0
 \end{array}\bigg) }
\overline{G_{k,r+1}} \ + \
{\small \bigg(\begin{array}{c}
n-r  \\
1
 \end{array}\bigg) }
\overline{G_{k,r}}
\ + \dots + \
{\small \bigg(\begin{array}{c}
n-2  \\
r-1
 \end{array}\bigg) }
\overline{G_{k,2}} \ \ .
$$
\end{proof}

Remark: In particular, Theorem 2 implies that the $M_{k,r}$
are Schur-convex. This fact can also be verified directly,
by checking (\ref{Sch}) with $\Phi = -M_{k,r}$.

Lastly, for completeness we mention without proof a result from \cite{k3}
indicating that there
is actually some  ``meaningful" property implied by the simultaneous assumptions
\begin{equation}
\label{F}
 F_{k,r}(x) \leq  F_{k,r}(y) \ \ \ \ {\rm for \ all \ integers}\ \  k,r\geq 1,
\end{equation}
although it is not the usual majorization relation. Namely, (\ref{F}) implies that
\begin{equation}
\label{L}
\sum_{i=1}^n \psi(x_i) \leq \sum_{i=1}^n \psi(y_i)
\end{equation}
for all $\psi : [0,\infty)\to [0,\infty)$ of the form
$ \psi(s) = \int_0^s \ \varphi(t) \frac{dt}{t}$
where $\varphi$ is concave nondecreasing, or equivalently, for all $\psi$ of
the form
$\psi(s) = \psi_\lambda(s)
:= \min(s,\lambda) + \lambda \log_+(s/\lambda), \ \lambda > 0$.
The proof is given in \cite[Theorem 9]{k3} but will not
be included in the present paper for lack of space. The converse implication is almost obtained
as well, except for a ``technical lemma"  which still requires proof \cite[Theorem 15, Conjecture 16]{k3}.
Thus it would appear that the simultaneous
inequalities (\ref{F}) are characterized by (\ref{L}).
%



\end{document}